%\magnification=\magstep1
\input amstex
\documentstyle{amsppt}
\topmatter
\NoBlackBoxes
\title The Largest Prime Dividing the Maximal Order \\
of an Element of $S_n$ \endtitle
\affil Department of Mathematics \\
       University of Georgia \\
       Athens, GA  30602 
\endaffil
\thanks
A portion of this research, including the computations, was done at
the Supercomputing Research Center.  The author would also like to
thank the referee for helpful suggestions.
\endthanks
\abstract   We define $g(n)$ to be the maximal order of an element of
the symmetric group on $n$ elements.  Results about the prime factorization of
$g(n)$ allow a reduction of the upper bound on the largest prime
divisor of $g(n)$ to $1.328\sqrt{n\log n}$.
\endabstract
\email grantham\@joe.math.uga.edu \endemail
\author Jon Grantham \endauthor
\define\landau{1}
\define\mass{2}
\define\mnr{3}
\define\nic{4}
\define\nica{5}
\define\nicb{6}
\define\rands{7}
\define\schoen{8}
\endtopmatter
\document

Let $S_n$ be the symmetric group on $n$ letters.  
\proclaim{Definition}
$g(n)=\operatorname{max}\,\{\operatorname{ord}(\sigma)\,|\,\sigma\in S_n\}.$
\endproclaim

The first work on $g(n)$ was done by Landau \cite{\landau}
in 1903.  He showed that $\log{g(n)} \sim \sqrt{n\log n}$ as
$n\to\infty$.
In 1984, 
Massias \cite{\mass} showed an upper bound for $\frac{\log g(n)}%
{\sqrt{n\log n}}$,
$$ \log g(n) \le a\sqrt{n\log n}\qquad a = 1.05313\dots\qquad n\ge 1,$$
with $a$ attained for $n = 1,319,166$.

Let $P(g(n))$ be the largest prime divisor of $g(n)$. In 1969,
Nicolas \cite{\nic}
proved that $P(g(n)) \sim \sqrt{n\log n}$ as $n\to\infty$.
In 1989, Massias, Nicolas, and Robin \cite{\mnr} showed that $P(g(n))\le
2.86\sqrt{n\log n}$, $n\ge 2$.
They conjectured that $\frac{P(g(n))}{\sqrt{n\log n}}$
achieves a maximum ($1.265\dots$) for $n\ge 5$ at $n=215$, with $P(g(215))=43$.
They note that improving this bound using the 
techniques of their proof would require ``very extensive computation,''
and even then would not be able to reduce the constant in the bound
below $2$.

Using different techniques, however, we can improve this result to the
following

\proclaim{Theorem}
For each integer $n\ge 5$, we have
$$P(g(n))\le 1.328 \sqrt{n\log n}.$$
\endproclaim

Our proof begins with the simple observation that
$g(n)=\operatorname{max}\,\{\operatorname{ord}(\sigma)\,|\,\sigma\in S_n'\},$
where $S_n'$ is the subset of $S_n$ consisting 
of elements that are the product of disjoint cycles of prime power length.

To see this, recall the fact that we can write any $\sigma\in S_n$ as the
product of disjoint cycles.  Then $\operatorname{ord}(\sigma)$ is the least
common multiple of the cycle lengths. Consider a cycle of 
length $ab$ with $(a,b)=1$, $a,b>1$. The product of a cycle of length
$a$ with one of length $b$ also has order $ab$ and is a permutation on
fewer elements.  Thus, given any element of $S_n$, we may find
another that has the same order and is a product of disjoint cycles of prime
power length.

\proclaim{Definition}
For each natural number $M$, let $\ell(M)=\sum\limits_{p^\alpha\|M}{p^\alpha}.$
\endproclaim
We observe that $\ell(M)$ is the shortest length of a permutation of
order $M$.
Thus, we can characterize $g(n)$ in terms of $\ell$ as follows:
$$g(n)=\text{max}\,\{M\,|\,\ell(M)\le n\}.$$
In particular, $\ell(g(n))\le n$.

Nicolas \cite{\nicb} describes an algorithm for computing $g(n)$.
Employing a variation of this algorithm, I computed exact values of $g(n)$ 
for $n\le 500,000$ on a Sun 4/390.  The accuracy of
the computation was checked by calculating values of $g(n)$ using the set
$G$ described in \cite{\mnr} and verifying that they matched
those in the computations.  Analysis of the computations confirmed that
for $5\le n\le 500,000$, $\frac{P(g(n))}{\sqrt{n\log{n}}}$ attains a
maximum at $n=215$.

\proclaim{Lemma 1 (Nicolas \cite{\nica})} 
Let $p$, $p'$, and $q$ be distinct primes, with $q \ge p+p'$.  If $q$
divides $g(n)$, then at least one of $p$ and $p'$ divides $g(n)$.
\endproclaim
\demo{Proof}
Suppose $p$ and $p'$ are primes not dividing $g(n)$.
Assume there is a prime $q\ge{p+p'}$ with $q\mid g(n)$.
Without loss of generality, $p<p'$.  Choose $k$ such that
$$p^k+p'\leq q \le p^{k+1}+p'-1.$$
Let $M=\frac{p^kp'g(n)}q$.  Since $q\mid g(n)$, $M$ is an integer.  Then
$$\ell(M)\le\ell(g(n))+(p^k+p'-q)\leq\ell(g(n))\leq n.$$
Thus, an element of order $M$ can be written as a permutation on $n$
letters.
Also,
$$ \multline
p^kp'-q\ge p^kp'-p^{k+1}-p'+1=p^k(p'-p)-p'+1 \\
\ge p(p'-p)-p'+1=(p-1)(p'-p-1)\ge0.
\endmultline$$
Therefore, $p^kp'>q$, so $M>g(n)$.  But $g(n)$ is the maximal order of a
permutation on $n$ letters.  Thus, we have a contradiction, and the lemma is
proven. %
\enddemo

Write $q=P(g(n))$.  We immediately get the following
\proclaim{Corollary}
At most one prime less than $\frac q2$ fails to divide $g(n)$.
\endproclaim

\proclaim{Lemma 2}
Suppose $0<\alpha<\beta<1$.  
If at least one prime in the interval $(\alpha q,\beta q)$
divides $g(n)$, then at most one prime in the interval
$(\sqrt\beta q, \frac{(1+\alpha)q}2)$ fails to divide $g(n)$.
\endproclaim

\demo{Proof}
If two primes in the interval $(\sqrt\beta q, \frac{(1+\alpha)q}2)$ fail
to divide $g(n)$, call them $p$ and $p'$.  Let $q'$ be a prime in the
interval $(\alpha q,\beta q)$ dividing $g(n)$.  Let
$M=\frac{pp'}{qq'}g(n)$.  Then
$$\ell(M)\le p+p'-q-q'+\ell(g(n))\le (1+\alpha)q-q-\alpha q+\ell(g(n))%
=\ell(g(n)).$$
But $pp'-qq'>(\sqrt\beta q)^2-q(\beta q)=0,$ so $M>g(n)$, giving a
contradiction.
\enddemo

\demo{Proof of Theorem}
By the computations, we may take $n>500\text{,}000$.  We may also
assume $q\ge 1.3\sqrt{500000\log 500000}>3329$.   Using the results of
Schoenfeld \cite{\schoen} for large $q$, and computations for small $q>3329$, we
see that there are always at least two primes in the
intervals  $(\alpha_iq,\beta_iq)$, with $\alpha_1=.2426$,
$\beta_1=.25$, $\alpha_2=.3746$, $\beta_2=.386$, 
$\alpha_3=.4632$, $\beta_3=.4723$,
$\alpha_4=.5248$, $\beta_4=.5352$,
$\alpha_5=.57$, $\beta_5=.5812$,
$\alpha_6=.6044$, $\beta_6=.6162$,
$\alpha_7=.6312$, $\beta_7=.6435$,
$\alpha_8=.6534$, $\beta_8=.6652$, and
$\alpha_9=.6714$, $\beta_9=.6834$.
By Lemma 1, at most one of the two or more primes in any of the first three intervals
fails to divide $g(n)$.
Applying Lemma 2, we get that at most one prime in each interval 
$(\sqrt{\beta_i}q,\frac{(1+\alpha_i)q}2)$ fails to divide $g(n)$, for
$i\le 3$.
This fact in turn implies that at most one prime in each interval $(\alpha_i
q, \beta_iq)$ fails to divide $g(n)$ for $4\le i\le 9$.  Applying Lemma 2
again, we see that at most one prime in each interval
$(\sqrt{\beta_i}q,\frac{(1+\alpha_i)q}2)$ fails to divide $g(n)$ for
$1\le i\le 9$.

We note that these intervals cover $(.5q,.8357q)$.  So at most ten
primes less than $.8357q$ fail to divide $g(n)$, and they can be at most
$\frac q2$, $\frac{(1+\alpha_1)q}2$,\dots, and $\frac{(1+\alpha_9)q}2$.

Therefore,
$$g(n)\ge \frac{q\prod\limits_{p\le .8357q}p}%
{\frac q2\prod {\frac{1+\alpha_i}2q}}.$$
Taking logarithms, we get
$$\log{g(n)}\ge \theta(.8357q)-\log{\dfrac 12}-\sum{\log\left(%
\frac{1+\alpha_i}2q\right)},$$
where $\theta$ is the Chebyshev function,
$\theta(x)=\sum\limits_{p\le x}{\log p}$.

For $q>3329$ the sum of the terms on the right is less than $.01338q$, so
$$\log{g(n)}\ge \theta(.8357q)-.01338q.$$
Using the estimates for $\theta(x)$ in \cite{\rands}, we get
$$\log{g(n)}\ge .79307q.$$
From \cite{\mnr}, $1.05314\sqrt{n\log n}\ge \log{g(n)}$, so
$$1.328\sqrt{n\log n} \ge \frac{1.05314}{.79307}\sqrt{n\log n} \ge q.$$

It is likely that further computation would be able to show that 
$\frac{P(g(n))}{\sqrt{n\log{n}}}$ attains a
maximum at $n=215$ for \bf{all} $n\ge 5$.

\enddemo

\enddocument